\documentclass[11pt]{article}

\usepackage{fullpage,amsmath,epsfig}
\usepackage{graphicx}
\usepackage{caption2}

\begin{document}

\newtheorem{theorem}{Theorem}[section]
\newtheorem{prop}[theorem]{Proposition}
\newtheorem{lemma}[theorem]{Lemma}
\newtheorem{remark}[theorem]{Remark}
\newtheorem{corollary}[theorem]{Corollary}
\newtheorem{observation}[theorem]{Observation}
\newtheorem{claim}[theorem]{Claim}
\newtheorem{conj}[theorem]{Conjecture}
\newtheorem{definition}[theorem]{Definition}
\def\square{\vrule height6pt width7pt depth1pt}
\def\endpf{\hfill\square\bigskip}

\title{Lattice Grids and Prisms are Antimagic}

\author{
\\
  Yongxi~Cheng
\\
\\
\small{Department of Computer Science,
  Tsinghua University, Beijing 100084,
  China}
  \\
[1mm] \small{cyx@mails.tsinghua.edu.cn}
}

\date{} 

\maketitle

\begin{abstract}
An \emph{antimagic labeling} of a finite undirected simple graph
with $m$ edges and $n$ vertices is a bijection from the set of edges
to the integers $1,\ldots,m$ such that all $n$ vertex sums are
pairwise distinct, where a vertex sum is the sum of labels of all
edges incident with the same vertex. A graph is called
\emph{antimagic} if it has an antimagic labeling. In 1990, Hartsfield and
Ringel conjectured that every connected graph, but $K_2$, is antimagic. In 2004, N. Alon et al
showed that this conjecture is true for $n$-vertex graphs with minimum degree
$\Omega(\log n)$. They also proved that complete partite graphs
(other than $K_2$) and $n$-vertex graphs with maximum degree at least $n-2$ are
antimagic. Recently, Wang showed that the toroidal grids (the Cartesian products of two
or more cycles) are antimagic. Two open problems left in
Wang's paper are about the antimagicness of lattice grid graphs and prism graphs, which are the Cartesian
products of two paths, and of a cycle and a path, respectively. In
this article, we prove that these two classes of
graphs are antimagic, by constructing such antimagic labelings.
\\[2mm]
\emph{Keywords:} Antimagic; Labeling; Lattice grid; Prism
\end{abstract}

\section{Introduction}
All graphs in this paper are finite, undirected and simple. In 1990,
Hartsfield and Ringel \cite{HaRi} introduced the concept of
\emph{antimagic} graph. An \emph{antimagic labeling} of a graph with $m$ edges and
$n$ vertices is a bijection from the set of edges to the integers
$1,\ldots,m$ such that all $n$ vertex sums are pairwise distinct,
where a vertex sum is the sum of labels of all edges incident with
that vertex. A graph is called antimagic if it has an antimagic
labeling. Hartsfield and Ringel showed that paths $P_n (n\geq 3)$,
cycles, wheels, and complete graphs $K_n (n\geq 3)$ are antimagic.
They conjectured that all trees except $K_2$ are antimagic.
Moreover, all connected graphs except $K_2$ are antimagic. These two
conjectures are unsettled. In 2004, Alon et al \cite{AKLRY} showed
that the latter conjecture is true for all graphs with $n$ vertices
and minimum degree $\Omega (\log n)$. They also proved that a graph
$G$ with $n\ (\geq 4)$ vertices and maximum degree $\Delta (G)\geq
n-2$ is antimagic, and all complete partite graphs except $K_2$ are
antimagic. In \cite{Wa}, Wang showed that the toroidal
grids (the Cartesian products of two cycles) are antimagic, the author
also proved that all Cartesian products of an antimagic $k$-regular graph ($k>1$)
and a cycle (consequently Cartesian products of more than two cycles) are antimagic.
Two open problems left in \cite{Wa} are
about the antimagicness of lattice grid graphs and prism graphs, which are the Cartesian
products of two paths, and of a cycle and a path, respectively.

In this paper, we prove that these two classes of graphs are
antimagic, by constructing such antimagic labelings.
In contrast to toroidal grids, lattices and prisms have less symmetry (more local
structures), we will incorporate new strategies in our labeling.
Our main results are the following two theorems, which are proved in
Section 3 and Section 4, respectively.

\begin{theorem}\label{lattice}
\textit{All lattice grid graphs $P_1[m+1]\times P_2[n+1]$ are
antimagic, for integers $m,n\geq 1$.}
\end{theorem}

\begin{theorem}\label{prism}
\textit{All prism graphs $C[m]\times P[n+1]$ are antimagic, for
integers $m\geq 3, n\geq 1$.}
\end{theorem}

For more results, open problems and conjectures on antimagic graphs
and various graph labeling problems, please see \cite{Ga, He}.

\section{Preliminaries}\label{prelim}
The \emph{Cartesian product} $G_1\times G_2$ of two graphs $G_1=(V_1, E_1)$
and $G_2=(V_2, E_2)$ is a graph with vertex set $V_1\times V_2$, and
$(u_1,u_2)$ is adjacent to $(v_1,v_2)$ in $G_1\times G_2$ if and
only if $u_1=v_1$ and $u_2v_2\in E_2$, or, $u_2=v_2$ and $u_1v_1\in
E_1$. The Cartesian product of two paths is a lattice grid graph,
and the Cartesian product of a path and a cycle is a prism grid
graph.

\begin{figure}[t]
\renewcommand{\captionlabelfont}{\bf}
\renewcommand{\captionlabeldelim}{.~}
\centering
\includegraphics[width=120mm]{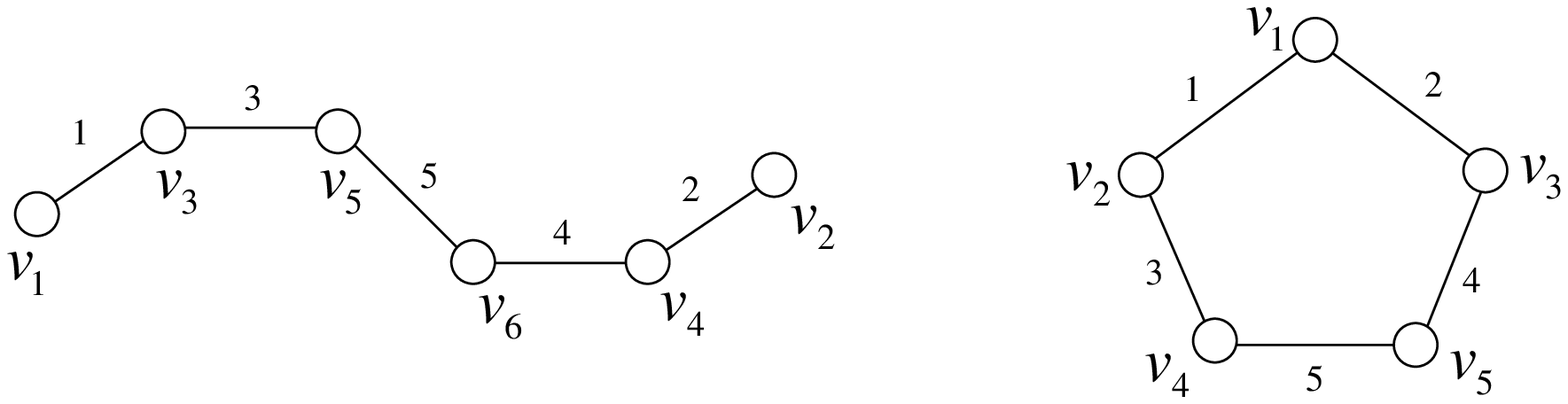}
\renewcommand{\figurename}{Fig.}
\caption{Antimagic labeling of $P[n+1]$ and $C[m]$, for $n=5$,
$m=5$} \label{fig:pathcycle}
\end{figure}

Before proving our main results, we first describe antimagic
labeling on paths and cycles respectively (see Figure \ref{fig:pathcycle}).
The labeling methods are the same as in \cite{Wa},
here we rephrase them for the sake of completeness.

\begin{lemma}\label{path}
\textit{All paths $P[m+1]$ are antimagic for integers $m\geq 2$.}
\end{lemma}

\noindent {\bf Proof:}\, Suppose the vertex set is $\{v_1,\ldots,
v_{m+1}\}$ and the edge set is arranged to be
$\{v_iv_{i+2}|i=1,\ldots,m-1\}\cup\{v_mv_{m+1}\}$. The following
labeling $f(v_iv_{i+2})=i$, for $1\leq i\leq m-1$, and
$f(v_mv_{m+1})=m$ is antimagic, since we have
$$
f^{+}(v_i) = \left\{
    \begin{array}{ll}
i & i=1,2;
\\
2i-2 & i=3,\ldots,m;
\\
2m-1 & i=m+1.
\end{array}
\right.
$$
Therefore,
$$
f^{+}(v_1)<f^{+}(v_2)< \ldots\ldots <f^{+}(v_{m+1})
$$\endpf

\begin{lemma}\label{cycle}
\textit{All cycles $C[m]$ are antimagic for integers $m\geq 3$.}
\end{lemma}

\noindent {\bf Proof:}\, Suppose the vertex set is $\{v_1,\ldots,
v_m\}$ and the edge set is arranged to be $\{v_1v_2\}\cup
\{v_iv_{i+2}|i=1,\ldots,m-2\}\cup\{v_{m-1}v_m\}$. The following
labeling $f(v_1v_2)=1$, $f(v_iv_{i+2})=i+1$, for $1\leq i\leq m-2$,
and $f(v_{m-1}v_m)=m$ is antimagic, since we have
$$
f^{+}(v_i) = \left\{
    \begin{array}{ll}
3 & i=1;
\\
2i & i=2,\ldots,m-1;
\\
2m-1 & i=m.
\end{array}
\right.
$$
Therefore,
$$
f^{+}(v_1)<f^{+}(v_2)< \ldots\ldots <f^{+}(v_m)
$$\endpf

\section{Proof of Theorem \ref{lattice}}

Let $f: E(P_1[m+1]\times P_2[n+1])\rightarrow
\{1,2,\ldots,2mn+m+n\}$ be an edge labeling of $P_1[m+1]\times
P_2[n+1]$, and denote the induced sum at vertex $(u,v)$ by
$f^{+}(u,v)=\sum f((u,v),(y,z))$ , where the sum runs over all
vertices $(y,z)$ adjacent to $(u,v)$ in $P_1[m+1]\times
P_2[n+1]$. To prove Theorem \ref{lattice}, first, we construct a
labeling that is antimagic on product graphs of two paths $P_1[m+1]$
and $P_2[n+1]$, for $n\geq m\geq 2$. Then, we give an antimagic
labeling of graphs $P_1[2]\times P_2[n+1]$, for $n\geq 1$.

\subsection{$P_1[m+1]\times P_2[n+1]$ is Antimagic, for $n\geq m\geq 2$}

Assume that $P_1[m+1]$ has edge set
$\{u_iu_{i+2}|i=1,\ldots,m-1\}\cup\{u_mu_{m+1}\}$, and $P_2[n+1]$
has edge set $\{v_iv_{i+1}|i=1,\ldots,n\}$. We will construct an antimagic labeling of
$P_1[m+1]\times P_2[n+1]$ for $n\geq m\geq 2$, which contains
two phases.
\\[2mm]
\noindent \emph{Phase 1:} For the $mn+m$ edges
contained in copies of $P_1[m+1]$ component (i.e., the edges
$((u_i,v_j), (u_{i+2},v_j))$ and $((u_m,v_j), (u_{m+1},v_j))$, for
$1\leq i\leq m-1, 1\leq j\leq n+1$), label them with even
numbers $2,4,\ldots,2mn+2m$ (notice $n\geq m$).

Specifically, first label the edges of $P_1[m+1]$ with
$U$ and $R$ such that $u_1u_3$ is labeled with $U$, and two edges
are labeled with different letters if they are incident to a same vertex. Obviously, there is one unique such
labeling. For each edge $u_iu_j\in E(P_1[m+1])$ labeled with $U$,
label the edges $((u_i,v_1),(u_j,v_1)), ((u_i,v_2),(u_j,v_2)),\ldots\ldots,((u_i,v_{n+1}),(u_j,v_{n+1}))$
in usual order; for each edge $u_iu_j\in E(P_1[m+1])$
 labeled with $R$, label the edges $((u_i,v_1),(u_j,v_1)),
((u_i,v_2),(u_j,v_2)),\ldots\ldots$, \\$((u_i,v_{n+1}),(u_j,v_{n+1}))$
in reversed order, and

\begin{center}
\noindent $2,\ 4,
\ldots\ldots\ldots\ldots\ldots\ldots\ldots,
2n+2$, \ \
(labels for $((u_1,v_i),(u_3,v_i)), i=1,2,\ldots,n+1$)
\\[2mm]
$2n+4,\ 2n+6,
\ldots\ldots\ldots\ldots,
4n+4$,\ \
(labels for $((u_2,v_i),(u_4,v_i)), i=1,2,\ldots,n+1$)
\\[2mm]
\ldots\ldots\ldots\ldots\ldots\ldots\ldots\ldots\ldots\ldots\ldots\ldots\ldots\ldots\ldots\ldots\ldots
\ldots\ldots\ldots\ldots\ldots\ldots\ldots\ldots\ldots
\\[2mm]
$2mn+2m-2n,\ldots, 2mn+2m$, \ \
(labels for $((u_m,v_i),(u_{m+1},v_i)), i=1,2,\ldots,n+1$)
\end{center}

\noindent \emph{Phase 2:} Denote by $A: a_1<a_2<\ldots<a_s$
the sequence of all odd numbers in $\{1,2,\ldots,2mn+m+n\}$, and
denote by $B: b_1<\ldots<b_t$ the sequence of all even numbers
in $\{2mn+2m+1,\ldots,2mn+m+n\}$, i.e., the even numbers that are
not used in Phase 1. Notice that $t\leq \frac
{1}{2}(2mn+m+n)-(mn+m)=\frac {1}{2}(n-m)$. We merge $A$ and $B$ into
a sequence $C:$ $a_1,a_2,\ldots,a_{s-t},b_1,a_{s-t+1},b_2,\ldots,b_t,a_s$ of $s+t$
terms ($s+t=mn+n$), and denote the sequence $C$ by $c_1,c_2,...,c_{mn+n}$, which are
the labels for the other $mn+n$ edges contained in copies of $P_2[n+1]$ component.

\begin{figure}[t]
\renewcommand{\captionlabelfont}{\bf}
\renewcommand{\captionlabeldelim}{.~}
\centering
\includegraphics[width=128mm]{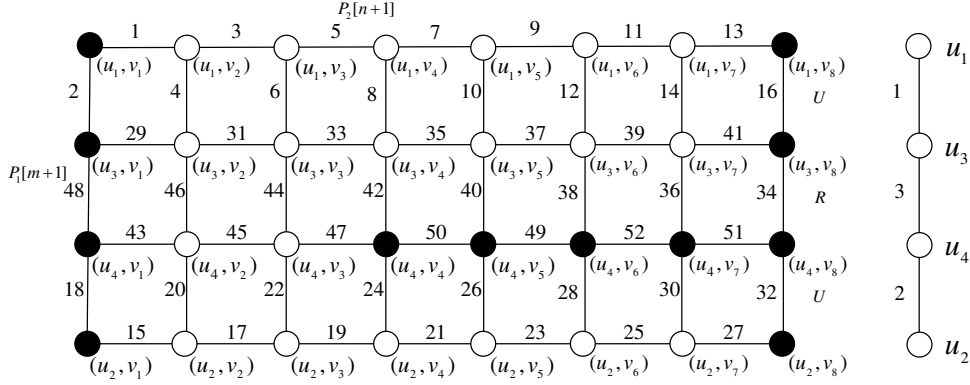}
\renewcommand{\figurename}{Fig.}
\caption{Antimagic labeling of $P_1[m+1]\times P_2[n+1]$, for $m=3,
n=7$} \label{fig:l48}
\end{figure}

For the $i$-th $P_2[n+1]$ component (with vertices $(u_i,v_1)$,
$(u_i,v_2)$,\ldots, $(u_i,v_{n+1})$), label its edges
in usual order according to the indices in the sequence $C$, $i=1,2,\ldots, m+1$,  and
\begin{center}
\noindent $c_1,\ c_2,
\ldots\ldots\ldots\ldots\ldots\ldots\ldots,
c_n$, \ \ (labels for the 1st $P_2[n+1]$ component)
\\[2mm]
$c_{n+1},\ c_{n+2},
\ldots\ldots\ldots\ldots\ldots,
c_{2n}$, \ \ (labels for the 2nd $P_2[n+1]$ component)
\\[2mm]
\ldots\ldots\ldots\ldots\ldots\ldots\ldots\ldots\ldots\ldots\ldots\ldots\ldots
\ldots\ldots\ldots\ldots\ldots\ldots\ldots\ldots\ldots
\\[2mm]
$c_{mn+1},c_{mn+2},
\ldots,c_{mn+n}$, \ \ (labels for the $(m+1)$-th $P_2[n+1]$ component)
\end{center}
Notice that $2t\leq n-m$, hence only the edges in the $(m+1)$-th $P_2[n+1]$ component
may be labeled with even numbers (see Figure~\ref{fig:l48}).

\

In what follows, we will show that the above labeling is antimagic.
In the product graph $P_1[m+1]\times P_2[n+1]$, at each vertex
$(u,v)$, the edges incident to this vertex can be partitioned into
two parts, one part is contained in a copy of $P_1[m+1]$ component,
and the other part is contained in a copy of $P_2[n+1]$ component.
Let $f^{+}_1(u,v)$ and $f^{+}_2(u,v)$ denote the sum at vertex
$(u,v)$ restricted to $P_1[m+1]$ component and $P_2[n+1]$ component
respectively, i.e., $f^{+}_1(u,v)=\sum f((u,v),(y,v))$, where the
sum runs over all vertices $y$ adjacent to $u$ in $P_1[m+1]$, and
$f^{+}_2(u,v)=\sum f((u,v),(u,z))$, where the sum runs over all
vertices $z$ adjacent to $v$ in $P_2[n+1]$. Therefore,
$f^{+}(u,v)=f^{+}_1(u,v)+f^{+}_2(u,v)$. The following two claims imply the antimagicness of the above labeling.

\begin{claim} \label{clai:lattice1}
\noindent For the above labeling of $P_1[m+1]\times P_2[n+1]$,
$n\geq m\geq 2$, we have
\begin{eqnarray*}
&&f^{+}(u_1,v_2)<f^{+}(u_1,v_3)<\ldots\ldots\ldots\ldots<f^{+}(u_1,v_n)<
\\
&&f^{+}(u_2,v_2)<f^{+}(u_2,v_3)<\ldots\ldots\ldots\ldots<f^{+}(u_2,v_n)<
\\
&&\;\;\;\;\;\;
\ldots\ldots\ldots\ldots\ldots\ldots\ldots\ldots\ldots\ldots\ldots\ldots\ldots\ldots\ldots
\\
&&f^{+}(u_m,v_2)<\,f^{+}(u_m,v_3)<\,\ldots\ldots\ldots<\,f^{+}(u_m,v_n)<
\\
&&f^{+}(u_{m+1},v_2)<\ldots<f^{+}(u_{m+1},v_{n-2t}),
\end{eqnarray*}
where $t$ $(\leq \frac {1}{2}(n-m))$ is the number of even numbers
in $\{2mn+2m+1,\ldots,2mn+m+n\}$. In addition, all the above sums
are even numbers.
\end{claim}
\noindent {\bf Proof:}\, Since $f^{+}_1(u_1,v_2)<f^{+}_1(u_1,v_3)<\ldots<f^{+}_1(u_1,v_n)$ and
$f^{+}_2(u_1,v_2)<f^{+}_2(u_1,v_3)<\ldots<f^{+}_2(u_1,v_n)$, we have $f^{+}(u_1,v_2)<f^{+}(u_1,v_3)<\ldots<f^{+}(u_1,v_n)$.
$f^{+}(u_1,v_n)<f^{+}(u_2,v_2)$ since
$f^{+}_1(u_1,v_n)<f^{+}_1(u_2,v_2)$ and
$f^{+}_2(u_1,v_n)<f^{+}_2(u_2,v_2)$. $f^{+}(u_2,v_2)<f^{+}(u_2,v_3)<
\ldots\ldots<f^{+}(u_2,v_n)$ since
$f^{+}_2(u_2,v_{i+1})-f^{+}_2(u_2,v_i)\geq 4$ and
$f^{+}_1(u_2,v_{i+1})-f^{+}_1(u_2,v_i)\geq -2$, it follows that $f^{+}(u_2,v_{i+1})-f^{+}(u_2,v_i)\geq 2$, for
$i=2,\ldots,n-1$. If $m=2$, $f^{+}_1(u_3,v_2)=f^{+}_1(u_3,v_n)>f^{+}_1(u_2,v_n)$;
if $m>2$, $f^{+}_1(u_3,v_2)>f((u_3,v_2),(u_j,v_2))>f((u_2,v_n),(u_4,v_n))=f^{+}_1(u_2,v_n)$, where $j=4$ or $5$.
Thus, in either case we have $f^{+}_1(u_2,v_n)<f^{+}_1(u_3,v_2)$. Clearly,
$f^{+}_2(u_2,v_n)<f^{+}_2(u_3,v_2)$. It follows that $f^{+}(u_2,v_n)<f^{+}(u_3,v_2)$.

For the vertices of degree $4$, clearly, $f^{+}_
1(u_i,v_2)=f^{+}_1(u_i,v_3)=\ldots\ldots\ldots\ldots=f^{+}_1(u_i,v_n)$ for $i=3,\ldots,m+1.$ Moreover, $f^{+}_
1(u_3,v_2)<f^{+}_1(u_4,v_2)<\ldots<f^{+}_1(u_{m+1},v_2)$ since $f((u_1,v_2),(u_3,v_2))<f((u_2,v_2),(u_4,v_2))<\ldots<
f((u_{m-1},v_2),(u_{m+1},v_2))<f((u_m,v_2),(u_{m+1},v_2))$. It follows that
\begin{eqnarray*}
&&f^{+}_
1(u_3,v_2)=f^{+}_1(u_3,v_3)=\ldots\ldots\ldots\ldots=f^{+}_1(u_3,v_n)<
\\
&&f^{+}_1(u_4,v_2)=f^{+}_1(u_4,v_3)=
\ldots\ldots\ldots\ldots=f^{+}_1(u_4,v_n)<
\\
&&\;\;\;\;\;\;\ldots\ldots\ldots\ldots\ldots\ldots\ldots\ldots\ldots\ldots\ldots\ldots\ldots\ldots\ldots
\\
&&f^{+}_1(u_m,v_2)=\,f^{+}_1(u_m,v_3)=\,
\ldots\ldots\ldots=\,f^{+}_1(u_m,v_n)<
\\
&&f^{+}_1(u_{m+1},v_2)=\ldots=f^{+}_1(u_{m+1},v_{n-2t}).
\end{eqnarray*}

On the other hand, since $c_1<c_2<\ldots<c_{mn+n-2t}$, we have that
\begin{eqnarray*}
&&f^{+}_2(u_3,v_2)<f^{+}_2(u_3,v_3)<\ldots\ldots\ldots\ldots<f^{+}_2(u_3,v_n)<
\\
&&f^{+}_2(u_4,v_2)<f^{+}_2(u_4,v_3)<
\ldots\ldots\ldots\ldots<f^{+}_2(u_4,v_n)<
\\
&&\;\;\;\;\;\;\ldots\ldots\ldots\ldots\ldots\ldots\ldots\ldots\ldots\ldots\ldots\ldots\ldots\ldots\ldots
\\
&&f^{+}_2(u_m,v_2)<\,f^{+}_2(u_m,v_3)<\,
\ldots\ldots\ldots<\,f^{+}_2(u_m,v_n)<
\\
&&f^{+}_2(u_{m+1},v_2)<\ldots<f^{+}_2(u_{m+1},v_{n-2t}).
\end{eqnarray*}

Therefore,
\begin{eqnarray*}
&&f^{+}(u_3,v_2)<f^{+}(u_3,v_3)<\ldots\ldots\ldots\ldots<f^{+}(u_3,v_n)<
\\
&&f^{+}(u_4,v_2)<f^{+}(u_4,v_3)<
\ldots\ldots\ldots\ldots<f^{+}(u_4,v_n)<
\\
&&\;\;\;\;\;\;\ldots\ldots\ldots\ldots\ldots\ldots\ldots\ldots\ldots\ldots\ldots\ldots\ldots\ldots\ldots
\\
&&f^{+}(u_m,v_2)<\,f^{+}(u_m,v_3)<\,
\ldots\ldots\ldots<\,f^{+}(u_m,v_n)<
\\
&&f^{+}(u_{m+1},v_2)<\ldots<f^{+}(u_{m+1},v_{n-2t}).
\end{eqnarray*}

All the above sums are even because each of them contains exactly two
odd labels. \endpf

\begin{claim}\label{clai:lattice2}
The remaining $2m+2+2t$ sums $f^{+}(u_1,v_1)$,
$f^{+}(u_1,v_{n+1})$, $f^{+}(u_2,v_1)$, $f^{+}(u_2,v_{n+1})$,\ldots,\\
$f^{+}(u_{m+1},v_1)$, $f^{+}(u_{m+1},v_{n+1})$, and
$f^{+}(u_{m+1},v_{n+1-2t})$, $f^{+}(u_{m+1},v_{n+2-2t})$,\ldots,
$f^{+}(u_{m+1},v_n)$  are pairwise distinct. In addition, they are
all odd numbers.
\end{claim}

\noindent {\bf Proof:}\, Let us first consider the $2m+2$ sums
$f^{+}(u_1,v_1)$, $f^{+}(u_1,v_{n+1})$, $f^{+}(u_2,v_1)$,
$f^{+}(u_2,v_{n+1})$,\ldots, $f^{+}(u_{m+1},v_1)$,
$f^{+}(u_{m+1},v_{n+1})$, there are two natural cases:
\\[2mm]
\emph{Case 1.} $m$ is odd. In this case $u_2u_4\in
E(P_1[m+1])$ is labeled with $U$, from the way we do the labeling, we
have $f^{+}_1(u_1,v_1)\leq f^{+}_1(u_1,v_{n+1})\leq
f^{+}_1(u_2,v_1)\leq f^{+}_1(u_2,v_{n+1})\leq\ldots\leq
f^{+}_1(u_{m+1},v_1)\leq f^{+}_1(u_{m+1},v_{n+1})$ and
$f^{+}_2(u_1,v_1)<f^{+}_2(u_1,v_{n+1})<f^{+}_2(u_2,v_1)<
f^{+}_2(u_2,v_{n+1})<\ldots<f^{+}_2(u_{m+1},v_1)<f^{+}_2(u_{m+1},v_{n+1}).$
Therefore, $
f^{+}(u_1,v_1)<f^{+}(u_1,v_{n+1})<f^{+}(u_2,v_1)<f^{+}(u_2,v_{n+1})<\ldots<f^{+}(u_{m+1},v_1)<f^{+}(u_{m+1},v_{n+1}).
$
\\[2mm]
\emph{Case 2.} $m$ is even. In this case $u_2u_j\in
E(P_1[m+1])$ is labeled with $R$ (where $j=3$ if $m=2$, $j=4$ if
$m>2$), the ordering of the $2m+2$ sums $f^{+}(u_1,v_1)$,
$f^{+}(u_1,v_{n+1})$, $f^{+}(u_2,v_1)$, $f^{+}(u_2,v_{n+1})$,\ldots,
$f^{+}(u_{m+1},v_1)$, $f^{+}(u_{m+1},v_{n+1})$ is the same as in
case 1, but between vertices $(u_2,v_1)$ and $(u_2,v_{n+1})$.
Specifically, we have $f^{+}_1(u_1,v_1)\leq f^{+}_1(u_1,v_{n+1})\leq
f^{+}_1(u_2,v_1), f^{+}_1(u_2,v_{n+1})\leq f^{+}_1(u_3,v_1)\leq
\ldots\leq f^{+}_1(u_{m+1},v_{n+1})$ and
$f^{+}_2(u_1,v_1)<f^{+}_2(u_1,v_{n+1})<f^{+}_2(u_2,v_1)<f^{+}_2(u_2,v_{n+1})
<\ldots<f^{+}_2(u_{m+1},v_1)<f^{+}_2(u_{m+1},v_{n+1}).$ Therefore,
$$
f^{+}(u_1,v_1)<f^{+}(u_1,v_{n+1})<f^{+}(u_2,v_1),f^{+}(u_2,v_{n+1})<\ldots<f^{+}(u_{m+1},v_1)<f^{+}(u_{m+1},v_{n+1}).
$$
Since
$f^{+}(u_2,v_1)=f^{+}_1(u_2,v_1)+f^{+}_2(u_2,v_1)=(4n+4)+(2n+1)=6n+5$,
and
$f^{+}(u_2,v_{n+1})=f^{+}_1(u_2,v_{n+1})+f^{+}_2(u_2,v_{n+1})=(2n+4)+(4n-1)=6n+3$,
it follows that
$f^{+}(u_1,v_1)<f^{+}(u_1,v_{n+1})<f^{+}(u_2,v_{n+1})<f^{+}(u_2,v_1)<\ldots<f^{+}(u_{m+1},v_1)<f^{+}(u_{m+1},v_{n+1}).
$\\

Thus, in any of the above two cases, the $2m+2$ sums
$f^{+}(u_1,v_1)$, $f^{+}(u_1,v_{n+1})$, $f^{+}(u_2,v_1)$,
$f^{+}(u_2,v_{n+1})$,\ldots, $f^{+}(u_{m+1},v_1)$,
$f^{+}(u_{m+1},v_{n+1})$ are pairwise distinct, and
$f^{+}(u_{m+1},v_{n+1})$ is the largest among them. For the other
$2t$ sums $f^{+}(u_{m+1},v_{n+1-2t})$,
$f^{+}(u_{m+1},v_{n+2-2t})$,\ldots, $f^{+}(u_{m+1},v_n)$, they are in strict increasing order
$f^{+}(u_{m+1},v_{n+1-2t})<f^{+}(u_{m+1},v_{n+2-2t})<\ldots<
f^{+}(u_{m+1},v_n)$, since
$f^{+}_1(u_{m+1},v_{n+1-2t})=f^{+}_1(u_{m+1},v_{n+2-2t})=\ldots=
f^{+}_1(u_{m+1},v_n)$ and
$f^{+}_2(u_{m+1},v_{n+1-2t})<f^{+}_2(u_{m+1},v_{n+2-2t})<\ldots<
f^{+}_2(u_{m+1},v_n)$.

At this point, the only remained issue is to
notice that $f^{+}(u_{m+1},v_{n+1-2t})>f^{+}(u_{m+1},v_{n+1})$,
since $f^{+}_1(u_{m+1},v_{n+1-2t})=f^{+}_1(u_{m+1},v_{n+1})$ and
$f^{+}_2(u_{m+1},v_{n+1-2t})=a_{s-t}+b_1\geq
(2mn+m+n-1-2t)+(2mn+2m+2)\geq
2mn+m+n-1-(n-m)+2mn+2m+2=4mn+4m+1>2mn+m+n\geq
a_s=f^{+}_2(u_{m+1},v_{n+1})$. Hence, the $2m+2t+2$ sums are
pairwise distinct. They are all odd numbers since each of them contains
exactly one odd label. \endpf

Combining Claim \ref{clai:lattice1} and Claim \ref{clai:lattice2},
we have proved that the above labeling of $P_1[m+1]\times P_2[n+1]$
is antimagic, for $n\geq m\geq 2$. Please see Figure \ref{fig:l48}
as an example of antimagic labeling of $P_1[m+1]\times P_2[n+1]$,
for $m=3, n=7$.

\subsection{$P_1[2]\times P_2[n+1]$ is Antimagic, for $n\geq 1$}

\begin{figure}[t]
\renewcommand{\captionlabelfont}{\bf}
\renewcommand{\captionlabeldelim}{.~}
\centering
\includegraphics[width=120mm]{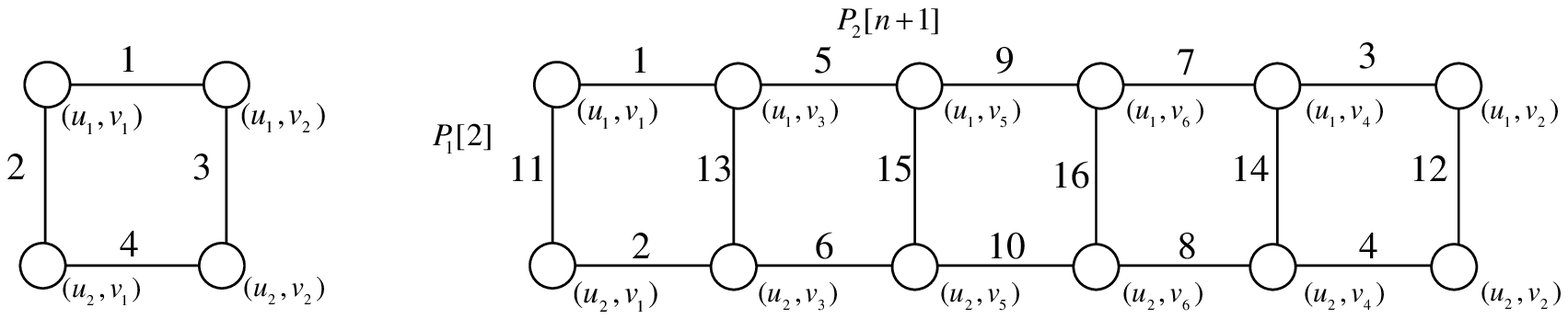}
\renewcommand{\figurename}{Fig.}
\caption{Antimagic labelings of $P_1[2]\times P_2[2]$ and $P_1[2]\times P_2[n+1]$, for $n=5$} \label{fig:l2n}
\end{figure}

Assume that $P_2[n+1]$ has edge set
$\{v_iv_{i+2}|i=1,\ldots,n-1\}\cup\{v_nv_{n+1}\}$. For $n=1$,
$P_1[2]\times P_2[2]$ is isomorphic to $C[4]$, hence by Lemma \ref{cycle}, it is antimagic. For $n>1$, label
$1,3,\ldots,2n-1$ to the edges $((u_1,v_1),(u_1,v_3))$,
$((u_1,v_2),(u_1,v_4))$,\ldots\ldots,
$((u_1,v_{n-1}),(u_1,v_{n+1}))$ ,$((u_1,v_n),(u_1,v_{n+1}))$, label
$2,4,\ldots,2n$ to the edges $((u_2,v_1),(u_2,v_3))$,
$((u_2,v_2),(u_2,v_4))$ ,\ldots\ldots,
$((u_2,v_{n-1}),(u_2,v_{n+1}))$, $((u_2,v_n),(u_2,v_{n+1}))$, and
label $2n+1,2n+2,\ldots,3n+1$ to $((u_1,v_1),(u_2,v_1))$,
$((u_1,v_2),(u_2,v_2))$, \ldots\ldots,
$((u_1,v_{n+1}),(u_2,v_{n+1}))$ (see Figure \ref{fig:l2n}).\\

We will show that the above labeling (for $n>1$) is antimagic.
Since the vertex sums restricted to $P_1[2]$ component satisfy that
$
f^{+}_1(u_1,v_1)=f^{+}_1(u_2,v_1)<f^{+}_1(u_1,v_2)=f^{+}_1(u_2,v_2)<\ldots
<f^{+}_1(u_1,v_{n+1})=f^{+}_1(u_2,v_{n+1})$ (`$=$' and `$<$'
alternate), and the vertex sums restricted to $P_2[n+1]$ component
are

$$
f^{+}_2(u_1, v_i) = \left\{
    \begin{array}{ll}
1 & i=1;
\\
3 & i=2;
\\
4i-6 & i=3,\ldots,n;
\\
4n-4 & i=n+1;
\end{array}
\right. ~~~~~~f^{+}_2(u_2, v_i) = \left\{
    \begin{array}{ll}
2 & i=1;
\\
4 & i=2;
\\
4i-4 & i=3,\ldots,n;
\\
4n-2 & i=n+1.
\end{array}
\right.
$$

It follows that $
f^{+}_2(u_1,v_1)<f^{+}_2(u_2,v_1)<f^{+}_2(u_1,v_2)<$\ldots$
<f^{+}_2(u_2,v_n)=f^{+}_2(u_1,v_{n+1})<f^{+}_2(u_2,v_{n+1})$ (there
is one equality). Therefore,
$f^{+}(u_1,v_1)<f^{+}(u_2,v_1)<f^{+}(u_1,v_2)<f^{+}(u_2,v_2)<\ldots
<f^{+}(u_1,v_{n+1})<f^{+}(u_2,v_{n+1})$, implying the antimagicness of the above labeling.\\

Combining the above two cases, we have proved Theorem \ref{lattice}.

\section{Proof of Theorem \ref{prism}}

Assume that in the product graph $C[m]\times P[n+1]$, $C[m]$ has
edge set $\{u_1u_2\}\cup
\{u_iu_{i+2}|i=1,\ldots,m-2\}\cup\{u_{m-1}u_m\}$, and $P[n+1]$ has
edge set $\{v_iv_{i+2}|i=1,\ldots,n-1\}\cup\{v_nv_{n+1}\}$. To prove Theorem \ref{prism}, first, we
construct a labeling that is antimagic on product graphs $C[m]\times
P[n+1]$ for $m\geq 3, n\geq 2$. Then, we give an antimagic labeling of
graphs $C[m]\times P[2]$ for $m\geq 3$.

\begin{lemma} \label{lemm:prism1}
\noindent $C[m]\times P[n+1]$ is antimagic for $m\geq 3, n\geq 2$.
\end{lemma}

\noindent {\bf Proof:}\, The antimagic labeling we will construct in
this case ($m\geq 3, n\geq 2$) is similar with the
labeling constructed in \cite{Wa} on toroidal grids, the
difference made here is to adapt the structure of prisms. The labeling contains two phases.
\\[2mm]
\noindent \emph{Phase 1:} Using the same way as in the antimagic labeling
of cycles in Lemma \ref{cycle}, label the edges on the $i$-th
$C[m]$ component (with vertices $(u_1,v_i)$,
$(u_2,v_i)$,\ldots, $(u_m,v_i)$), for $i=1,2,\ldots, n+1$, and
\begin{center}
\noindent $1,\ 2,
\ldots\ldots\ldots\ldots\ldots\ldots\ldots\ldots\ldots\ldots, m$, \
\  (labels for the 1st $C[m]$ component)
\\[2mm]
$m+1,\ m+2, \ldots\ldots\ldots\ldots\ldots\ldots\ldots, 2m$, \ \
(labels for the 2nd $C[m]$ component)
\\[2mm]
\ldots\ldots\ldots\ldots\ldots\ldots\ldots\ldots\ldots\ldots\ldots\ldots\ldots\ldots\ldots
\ldots\ldots\ldots\ldots\ldots\ldots\ldots\ldots\ldots
\\[2mm]
$mn+1,\ mn+2, \ldots\ldots\ldots mn+m$. \ \ (labels for the
$(n+1)$-th $C[m]$ component)
\end{center}

\begin{figure}[t]
\renewcommand{\captionlabelfont}{\bf}
\renewcommand{\captionlabeldelim}{.~}
\centering
\includegraphics[width=100mm]{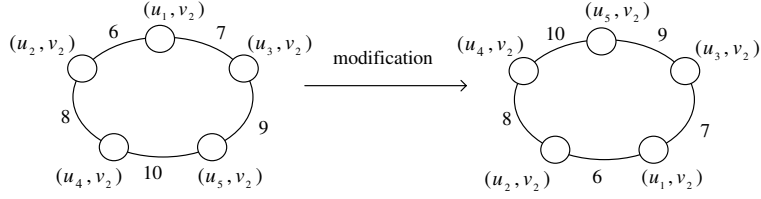}
\renewcommand{\figurename}{Fig.}
\caption{Modification on the 2nd $C[m]$ component in case $n$ is
even, for $m=5$} \label{fig:modi}
\end{figure}

\noindent \emph{Phase 2:} Similarly, label the edges of $P[n+1]$ with $U$ and $R$ such that
$v_1v_3$ is labeled with $U$, and two edges are labeled with
different letters if they are incident to a same vertex. For each edge
$v_iv_j\in E(P[n+1])$ labeled with $U$, the edges $((u_1,v_i),(u_1,v_j))$,
$((u_2,v_i),(u_2,v_j))$,\ldots\ldots,$((u_m,v_i),(u_m,v_j))$ will be
labeled in usual order; for each edge $v_iv_j\in E(P[n+1])$ labeled with $R$,
the edges $((u_1,v_i),(u_1,v_j))$, $((u_2,v_i),(u_2,v_j))$,\ldots\ldots,$((u_m,v_i),(u_m,v_j))$ will be
labeled in reversed order, and
\begin{center}
$mn+m+1,\ mn+m+2, \ldots,
mn+2m$, \ \ (labels for $((u_i,v_1),(u_i,v_3)), i=1,2,\ldots,m$)
\\[2mm]
$mn+2m+1,mn+2m+2, \ldots,
mn+3m$, \ (labels for $((u_i,v_2),(u_i,v_4)), i=1,2,\ldots,m$)
\\[2mm]
\ldots\ldots\ldots\ldots\ldots\ldots\ldots\ldots\ldots\ldots\ldots\ldots\ldots\ldots\ldots\ldots\ldots
\ldots\ldots\ldots\ldots\ldots\ldots\ldots\ldots\ldots
\\[2mm]
$2mn+1,2mn+2,
\ldots\ldots\ldots,
2mn+m$, \ \ (labels for $((u_i,v_n),(u_i,v_{n+1})), i=1,2,\ldots,m$)
\end{center}
If $v_2v_j\in E(P[n+1])$ ($j=3$ if $n=2$, $j=4$ if $n>2$) is
labeled with $R$ (i.e., when $n$ is even), we will take a
\emph{modification} process on the 2nd $C[m]$ component (with vertices $(u_1,v_2)$,
$(u_2,v_2)$,\ldots, $(u_m,v_2)$), which goes as follows. For each $u_iu_j\in E(C[m])$, the edge
$((u_i,v_2),(u_j,v_2))$ will be
relabeled with $(3m+1)-l_0(i,j)$, where $l_0(i,j)$ is the original
label assigned to $((u_i,v_2),(u_j,v_2))$ in Phase 1 (i.e., we `reverse' the labeling on the 2nd $C[m]$ component,
whose edges will still be labeled with the
same set of numbers $\{m+1,m+2,\ldots, 2m\}$).
Then, we rename each vertex $(u_i,v_2)$ as $(u_{m+1-i},v_2)$, for
$i=1,2,\ldots,m$ (see Figure \ref{fig:modi}).\\

Let $f^{+}_1(u,v)$ and $f^{+}_2(u,v)$ be the
vertex sum at $(u,v)\in V(C[m]\times P[n+1])$
restricted to $C[m]$ component and $P[n+1]$
component, respectively. Then,
$f^{+}(u,v)=f^{+}_1(u,v)+f^{+}_2(u,v)$ is the vertex sum at $(u,v)$.
It is easy to see that, for the above labeling, independent of the parity of $n$ (i.e.,
no matter whether there is a modification
process or not), the ordering
$f^{+}_1(u_1,v_2)<f^{+}_1(u_2,v_2)<\ldots\ldots<f^{+}_1(u_m,v_2)$
and
$f^{+}_2(u_1,v_2)<f^{+}_2(u_2,v_2)<\ldots\ldots<f^{+}_2(u_m,v_2)$
will hold.\\

Using similar arguments, it is straightforward to prove that for the above labeling we have
\begin{center}
$f^{+}_1(u_1,v_1)<f^{+}_1(u_2,v_1)<\ldots\ldots\ldots\ldots<f^{+}_1(u_m,v_1)<$
\\[2mm]
$f^{+}_1(u_1,v_2)<f^{+}_1(u_2,v_2)<\ldots\ldots\ldots\ldots<f^{+}_1(u_m,v_2)<$
\\
\ldots\ldots\ldots\ldots\ldots\ldots\ldots\ldots\ldots\ldots\ldots\ldots\ldots\ldots\ldots
\\
$f^{+}_1(u_1,v_{n+1})<f^{+}_1(u_2,v_{n+1})<\ldots\ldots<f^{+}_1(u_m,v_{n+1}),$
\end{center}

and

\begin{center}
$f^{+}_2(u_1,v_1)\leq f^{+}_2(u_2,v_1)\leq
\ldots\ldots\ldots\ldots\leq f^{+}_2(u_m,v_1)\leq $
\\[2mm]
$f^{+}_2(u_1,v_2)\leq f^{+}_2(u_2,v_2)\leq
\ldots\ldots\ldots\ldots\leq f^{+}_2(u_m,v_2)\leq $
\\
\ldots\ldots\ldots\ldots\ldots\ldots\ldots\ldots\ldots\ldots\ldots\ldots\ldots\ldots\ldots
\\
$f^{+}_2(u_1,v_{n+1})\leq f^{+}_2(u_2,v_{n+1})\leq \ldots\ldots\leq
f^{+}_2(u_m,v_{n+1}).$
\end{center}

Therefore,

\begin{center}
$f^{+}(u_1,v_1)<f^{+}(u_2,v_1)<\ldots\ldots\ldots\ldots<f^{+}(u_m,v_1)<$
\\[2mm]
$f^{+}(u_1,v_2)<f^{+}(u_2,v_2)<\ldots\ldots\ldots\ldots<f^{+}(u_m,v_2)<$
\\
\ldots\ldots\ldots\ldots\ldots\ldots\ldots\ldots\ldots\ldots\ldots\ldots\ldots\ldots\ldots
\\
$f^{+}(u_1,v_{n+1})<f^{+}(u_2,v_{n+1})<\ldots\ldots<f^{+}(u_m,v_{n+1}),$
\end{center}

\begin{figure}[t]
\renewcommand{\captionlabelfont}{\bf}
\renewcommand{\captionlabeldelim}{.~}
\centering
\includegraphics[width=100mm]{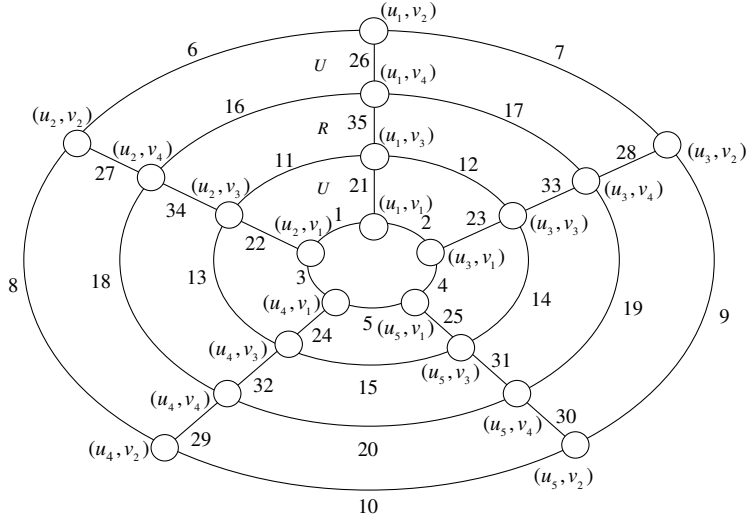}
\renewcommand{\figurename}{Fig.}
\caption{Antimagic labeling of $C[m]\times P[n+1]$, for $m=5$,
$n=3$} \label{fig:c5p4}
\end{figure}

which implies that the above labeling is antimagic. Please see
Figure \ref{fig:c5p4} as an example of antimagic labeling of
$C[m]\times P[n+1]$, for $m=5$, $n=3$.
\endpf

\begin{lemma} \label{lemm:prism2}
\noindent $C[m]\times P[2]$ is antimagic for $m\geq 3$.
\end{lemma}

\begin{figure}[t]
\renewcommand{\captionlabelfont}{\bf}
\renewcommand{\captionlabeldelim}{.~}
\centering
\includegraphics[width=60mm]{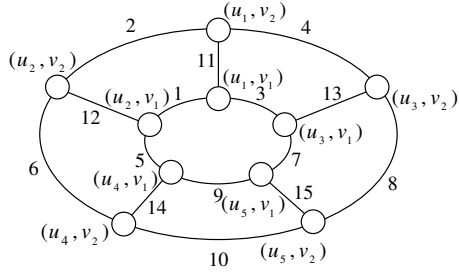}
\renewcommand{\figurename}{Fig.}
\caption{Antimagic labeling of $C[m]\times P[2]$, for $m=5$}
\label{fig:c5p2}
\end{figure}

\noindent {\bf Proof:}\, Assume that $C[m]$ has edge set
$\{u_1u_2\}\cup \{u_iu_{i+2}|i=1,\ldots,m-2\}\cup\{u_{m-1}u_m\}$.
Label $1,3,\ldots,2m-1$ to the edges $((u_1,v_1),(u_2,v_1))$,
$((u_1,v_1),(u_3,v_1))$,\ldots\ldots, $((u_{m-2},v_1),(u_m,v_1))$
,$((u_{m-1},v_1),(u_m,v_1))$, label $2,4,\ldots,2m$ to the edges
$((u_1,v_2),(u_2,v_2))$, $((u_1,v_2),(u_3,v_2))$,\ldots\ldots,
$((u_{m-2},v_2),(u_m,v_2))$ ,$((u_{m-1},v_2),(u_m,v_2))$, and label
$2m+1,2m+2,\ldots,3m$ to the edges \\
$((u_1,v_1),(u_1,v_2))$, $((u_2,v_1),(u_2,v_2))$, \ldots\ldots,
$((u_m,v_1),(u_m,v_2))$ (see Figure \ref{fig:c5p2}).\\

We will show that the above labeling ($m\geq 3$) is
antimagic. Since the vertex sums restricted to $C[m]$ component are

$$
f^{+}_1(u_i, v_1) = \left\{
    \begin{array}{ll}
4 & i=1;
\\
4i-2 & i=2,\ldots,m-1;
\\
4m-4 & i=m;
\end{array}
\right. ~~~~~f^{+}_1(u_i, v_2) = \left\{
    \begin{array}{ll}
6 & i=1;
\\
4i & i=2,\ldots,m-1;
\\
4m-2 & i=m.
\end{array}
\right.
$$
It follows that $
f^{+}_1(u_1,v_1)<f^{+}_1(u_1,v_2)=f^{+}_1(u_2,v_1)<\ldots
<f^{+}_1(u_{m-1},v_2)=f^{+}_1(u_m,v_1)<f^{+}_1(u_m,v_2)$ (there are
two equalities). In addition, $f^{+}_2(u_1,v_1)=f^{+}_2(u_1,v_2)<f^{+}_2(u_2,v_1)=f^{+}_2(u_2,v_2)<\ldots
<f^{+}_2(u_m,v_1)=f^{+}_2(u_m,v_2)$ (`$=$' and `$<$' alternate).
Therefore,
$f^{+}(u_1,v_1)<f^{+}(u_1,v_2)<f^{+}(u_2,v_1)<f^{+}(u_2,v_2)<\ldots
<f^{+}(u_m,v_1)<f^{+}(u_m,v_2)$, implying the antimagicness of the above labeling.
\endpf

Combining Lemma \ref{lemm:prism1} and Lemma \ref{lemm:prism2}, we
have proved Theorem \ref{prism}.

\end{document}